\documentclass[11pt,a4paper,twoside]{article}
\usepackage[UKenglish]{babel}
\usepackage[T1]{fontenc}
\usepackage[utf8]{inputenc}
\usepackage{latexsym,amsfonts,amsmath,amsthm,amssymb,mathrsfs}
\usepackage{fullpage}
\usepackage{xcolor,bm,bbm}
\usepackage{paralist}
\usepackage{fourier}
\usepackage{mathtools}
\usepackage{tikz}
\usetikzlibrary{patterns}
\mathtoolsset{showonlyrefs}

\newtheorem{thm}{Theorem}[section]

\newtheorem{proposition}[thm]{Proposition}
\newtheorem{rem}[thm]{Remark}

\newtheorem{lemma}[thm]{Lemma}

\newtheorem{thmalpha}{Theorem}

\allowdisplaybreaks
\newcommand{\R}{\mathbb R}
\newcommand{\N}{\mathbb N}
\newcommand{\M}{\mathbb M}
\newcommand{\Sp}{\mathcal S}

\newcommand{\rank}{\operatorname{rank}}
\newcommand{\Tr}{\operatorname{Tr}}
\newcommand{\codim}{\operatorname{codim}}
\newcommand{\HS}{\mathrm{HS}}
\newcommand{\embed}[2]{\Sp_{#1}^N\hookrightarrow\Sp_{#2}^N}

\begin{document}

\title{Gelfand numbers of Schatten class embeddings\\ in the quasi-Banach range}

\author{Joscha Prochno\footnote{University of Passau, Faculty of Computer Science and Mathematics, Dr.-Hans-Kapfinger-Stra{\ss}e 30, 94032 Passau, Germany. Email: joscha.prochno@uni-passau.de}}

\date{}
\maketitle

\begin{abstract}
We prove Conjecture~1.4 of Prochno and Strzelecki [J. Approx. Theory 277 (2022), 105736] on Gelfand numbers of embeddings between finite-dimensional Schatten classes. For $0<p\le1$ and indices bounded away from the ambient dimension by a fixed proportion, we obtain the conjectured lower bound, uniformly in the target exponent $2\le q\le\infty$. Together with the known upper estimate obtained by Hinrichs, Prochno and Vyb\'iral [Math. Ann., 380(3–4):1563–1593, 2021], this gives the exact asymptotic order in the previously unresolved range. As a consequence, we obtain a complete all-index description by combining our result with the estimates of Hinrichs, Prochno and Vyb\'iral and of Prochno and Strzelecki. The proof is deterministic and is based on Hilbert--Schmidt projections on rectangular matrix spaces.
\bigskip

\par\smallskip

\end{abstract}

\section{Introduction \& Main Results}\label{sec:introduction} 

The Schatten classes form one of the classical families of operator ideals in functional analysis. For a Hilbert space $H$ and $0<p<\infty$, the Schatten $p$-class $\Sp_p$ consists of all compact operators on $H$ whose sequence of singular values belongs to $\ell_p$. We write $\Sp_\infty$ for the compact operators equipped with the operator norm. The cases $p=1$ and $p=2$ give the trace-class, or nuclear, operators and the Hilbert--Schmidt operators, respectively. The origins of the theory lie in Schatten's work on cross-spaces of linear transformations \cite{S1946} and his subsequent papers with von Neumann \cite{SvN1946_II,SvN1946_III}. These investigations formed part of the development of nuclear operators and tensor products, with related contributions beyond the Hilbert space setting by Ruston \cite{R1951} and Grothendieck \cite{G233}. Schatten's monograph \cite{Schatten1960} subsequently provided a systematic account of norm ideals in terms of symmetric gauge functions and unitarily invariant crossnorms. As he observed there, the classical study of individual compact operators had only comparatively recently been complemented by an interest in the structure of spaces of such operators. 

This shift from operators to operator spaces naturally leads to questions about their geometry. For finite $p$, the Schatten class $\Sp_p$ is a non-commutative counterpart of $\ell_p$, sharing its nested scale $\Sp_p\subset\Sp_q$ for $p<q$, uniform convexity for $1<p<\infty$, and, in this range, duality with the conjugate exponent under the trace pairing; see \cite{Schatten1960,TJ1974}. Nevertheless, the geometry of Schatten classes is not determined by the corresponding sequence-space results. Building on the results of Kwapie\'n and Pe{\l}czy\'nski \cite{KP1970}, Gordon and Lewis \cite{GL1974} proved that, on an infinite-dimensional Hilbert space, the Banach spaces $\Sp_p$, $1\le p\le\infty$, fail to have local unconditional structure
whenever $p\ne2$. In finite dimensions, the geometry of Schatten unit balls has been investigated through isotropic constants \cite{KMP1998}, concentration and thin-shell phenomena \cite{GP2007,RV2016}, asymptotic volumes and volume ratios \cite{KPT2020_volume_ratio}, and volumes of intersections \cite{KPT2020_intersections}. A complementary probabilistic perspective is provided by Sanov-type large deviation principles for empirical eigenvalue and singular-value distributions of random matrices \cite{KPT2019_sanov}.
\newline
For approximation theory, a particularly relevant aspect of this geometry is the behaviour of the identity mappings $\Sp_p^N\hookrightarrow\Sp_q^N$ on subspaces of prescribed codimension. This is quantified by their Gelfand numbers. Geometrically, these numbers measure how small a section of the domain unit ball can be made in the target norm. Their determination therefore requires understanding the singular-value structure of matrices in subspaces of prescribed codimension and how it governs the comparison of different Schatten quasi-norms. This brings the distinction between commutative and non-commutative geometry directly into the problem: sharp estimates for the embeddings $\ell_p^N\hookrightarrow\ell_q^N$ do not, by themselves, settle the corresponding questions for matrices.
\newline
The geometry of low-codimensional matrix subspaces that underlies the study of Gelfand numbers also plays a central role in low-rank matrix recovery and completion, where Schatten norms and quasi-norms play a role analogous to that of sequence norms and quasi-norms in sparse recovery; see \cite{CR2009,CDK2015,FR2013,KS2018,RT2011}. In particular, the range $0<p<1$ is relevant to non-convex recovery by Schatten quasi-norm minimization. The connection with Gelfand numbers is more than a formal analogy: Ch\'avez-Dom\'inguez and Kutzarova \cite{CDK2015} characterized the stability of such recovery procedures in terms of Gelfand widths of suitable Schatten class embeddings. Sharp estimates for these quantities thus describe both the geometry of matrix subspaces and limitations on stable recovery from linear measurements.
\newline
The study of general approximation quantities for finite-dimensional sequence spaces has a long history. For Gelfand numbers, fundamental contributions are due to Gluskin \cite{G83} and Garnaev and Gluskin \cite{GG1984}, whereas the entropy numbers of embeddings between finite-dimensional $\ell_p$-spaces were studied by Sch\"utt \cite{Schuett1984} and subsequently in \cite{EdmundsTriebel1996,Triebel1997,GuedonLitvak2000,Kuehn2001}. In the non-commutative setting, related approximation quantities for embeddings between Schatten classes had already been investigated by Carl and Defant \cite{CD1997} and Hinrichs and Michels \cite{HM2005}, in addition to the work of Ch\'avez-Dom\'inguez and Kutzarova \cite{CDK2015} discussed above. Hinrichs, Prochno and Vyb\'iral \cite{HPV2017} subsequently determined the asymptotic behaviour of the entropy numbers of $\Sp_p^N\hookrightarrow\Sp_q^N$ for all $0<p,q\le\infty$, up to constants depending only on $p$ and $q$. However, the corresponding description of Gelfand numbers proved less complete. In \cite{HPV2021}, the authors obtained asymptotically sharp estimates in a wide range of parameters, extending classical sequence-space results and complementing the earlier bounds. Prochno and Strzelecki \cite{PS2022} subsequently established matching bounds for Gelfand numbers in the intermediate codimension range for $0<p\le2\le q\le\infty$ and obtained new estimates for approximation and Kolmogorov numbers of Schatten class embeddings. Still, two gaps remained in the description of Gelfand numbers in \cite{PS2022}. One concerns embeddings from $\Sp_p^N$, with $0<p<1$, into $\Sp_q^N$, with $2<q\le\infty$: sharp upper bounds were available, but a matching lower estimate throughout the range
of indices bounded away from $N^2$ by a fixed proportion had not been established. The purpose of the present paper is to prove the lower estimate proposed in \cite[Conjecture~1.4]{PS2022}. Together with the previously known bounds, including the boundary cases
$p=1$ and $q=2$, this completes the asymptotic description of the Gelfand numbers for all indices in the range $0<p\le1$ and $2\le q\le\infty$.

\subsection{Presentation of the main results} 

For $0<p\le\infty$ and $N\in\N$, let $\Sp_p^N$ denote the space of all real $N\times N$ matrices equipped with the Schatten $p$-norm (or quasi-norm when $0<p<1$), defined by
\[
    \|A\|_p
    :=
    \begin{cases}
        \displaystyle
        \left(\sum_{j=1}^{N}\sigma_j(A)^p\right)^{1/p}
        &: \, 0<p<\infty,\\[2mm]
        \sigma_1(A),
        &:\,  p=\infty,
    \end{cases}
\]
where $\sigma_1(A)\ge\cdots\ge\sigma_N(A)\ge0$ are the singular values of $A$, that is, the eigenvalues of $(A^{\mathsf T}A)^{1/2}$, counted with multiplicity. In particular, the Schatten $\infty$-norm coincides with the
operator norm of $A:\ell_2^N\to\ell_2^N$. We write $\Sp_p^N\hookrightarrow\Sp_q^N$ for the natural identity mapping between these spaces. Given a bounded linear operator $T:X\to Y$ between quasi-Banach spaces, its $n$-th Gelfand number is defined by
\begin{equation}\label{eq:def-gelfand}
 c_n(T)=\inf\bigl\{\bigl\|\left.T\right|_F\bigr\|:F\subset X\text{ closed and linear},\ \codim F<n\bigr\}.
\end{equation}
Equivalently, the Gelfand numbers are given by (here and below, all infima over $F$ are taken over linear subspaces, which are automatically closed in the finite-dimensional setting)
\[
 c_n(\embed pq)
 =\inf_{\substack{F\subset\R^{N\times N}\\\codim F<n}} 
 \sup_{A\in F\setminus\{0\}}\frac{\|A\|_q}{\|A\|_p},
 \qquad 1\le n\le N^2.
\]
Throughout, $1/\infty:=0$ and constants implicit in $\lesssim$ and $\asymp$ are independent of $N$ and $n$. Dependence on fixed parameters is indicated by subscripts.
\medskip

As already mentioned, the purpose of this paper is to determine the Gelfand numbers of $\Sp_p^N\hookrightarrow\Sp_q^N$ in the range $0<p\le1$ and $2\le q\le\infty$. Hinrichs, Prochno and Vyb\'iral \cite[Proposition~4.5]{HPV2021} obtained upper estimates for all indices in this range. Prochno and Strzelecki subsequently proved matching lower estimates in the intermediate range of codimensions \cite[Proposition~3.1]{PS2022}. For small and moderately large codimensions, the missing estimate was formulated in \cite[Conjecture~1.4]{PS2022}, which asserts that
\begin{equation}\label{eq:intro-conjecture}
 c_n(\embed pq)
 \gtrsim_{p,q}
 \min\left\{1,\frac Nn\right\}^{1/p-1/2},
 \qquad 1\le n\le(1-c_0)N^2;
\end{equation}
here $0<p\le1$, $2\le q\le\infty$, and $c_0\in(0,1)$ is the fixed universal constant denoted by $c$ in \cite{PS2022}. The endpoint $p=1$ was already covered by \cite[Proposition~3.2]{PS2022}, and the case $q=2$ was also known; see \cite[Proposition~3.3]{PS2022} and the references therein. Thus, the missing part of \eqref{eq:intro-conjecture} concerns $0<p<1$ and $q>2$.

Our main result establishes the conjectured lower estimate for any fixed proportion of the ambient dimension. We include the known matching upper estimate to emphasise that both bounds can be chosen uniformly in the target space index $q$.

\begin{thmalpha}\label{thm:main}
Let $0<p\le1$ and $\eta\in(0,1)$. There exist constants $c_{p,\eta}>0$ and $C_p>0$ such that, for every $2\le q\le\infty$, every $N\in\N$, and every integer $1\le n\le(1-\eta)N^2$,
\begin{equation}\label{eq:main}
 c_{p,\eta}\min\left\{1,\frac Nn\right\}^{1/p-1/2}
 \le c_n(\embed pq)
 \le C_p\min\left\{1,\frac Nn\right\}^{1/p-1/2}.
\end{equation}
\end{thmalpha} 

The lower estimate in Theorem~\ref{thm:main} is the contribution of this paper. Taking $\eta=c_0$ proves \cite[Conjecture~1.4]{PS2022}. The upper estimate is known and follows from the corresponding bound with target space $\Sp_2^N$ (see \cite[Proposition~3.3, with $q=2$]{PS2022}). In particular, its constant depends only on $p$, whereas the lower-bound constant depends only on $p$ and $\eta$. An explicit choice of $c_{p,\eta}$ is given in \eqref{eq:explicit-constant} below. 

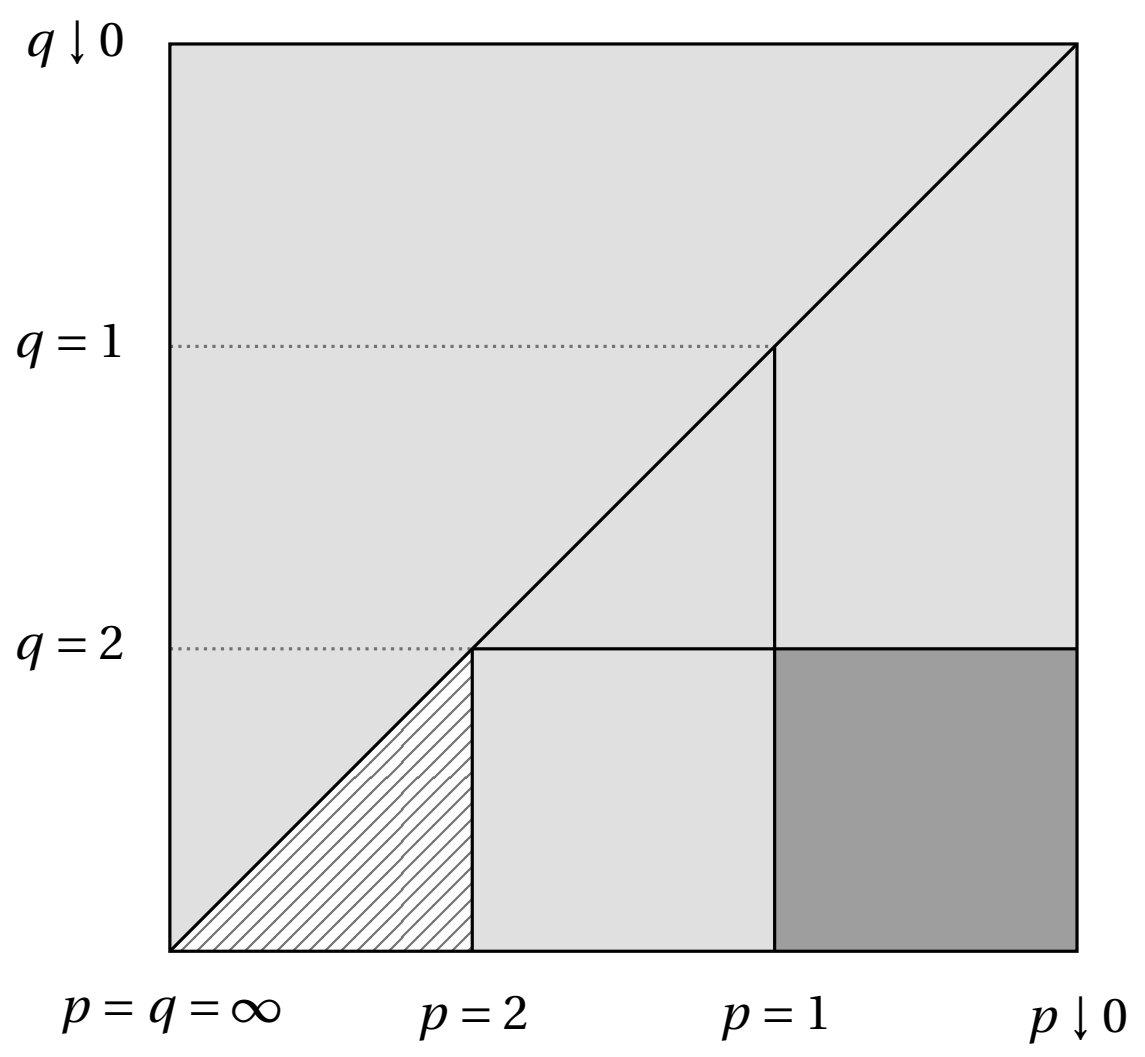
\begin{figure}[t] 
\centering
\begin{tikzpicture}[
    x=2cm, y=2cm,
    line width=0.6pt,
    font=\small
]

    \fill[black!12] (0,0) rectangle (3,3);

    \fill[black!38] (2,0) rectangle (3,1);

    \fill[white] (0,0) -- (1,0) -- (1,1) -- cycle;
    \fill[pattern=north east lines, pattern color=black!55]
        (0,0) -- (1,0) -- (1,1) -- cycle;

    \draw[densely dotted, black!55] (0,1) -- (1,1);
    \draw[densely dotted, black!55] (0,2) -- (2,2);

    \draw (0,0) rectangle (3,3);
    \draw (0,0) -- (3,3);
    \draw (1,0) -- (1,1);
    \draw (2,0) -- (2,2);
    \draw (1,1) -- (3,1);

    \node[below=5pt] at (0,0) {$p=q=\infty$};
    \node[below=5pt] at (1,0) {$p=2$};
    \node[below=5pt] at (2,0) {$p=1$};
    \node[below=5pt] at (3,0) {$p\downarrow0$};

    \node[left=5pt] at (0,1) {$q=2$};
    \node[left=5pt] at (0,2) {$q=1$};
    \node[left=5pt] at (0,3) {$q\downarrow0$};
\end{tikzpicture}

\caption{Schematic overview of the Gelfand numbers of $\Sp_p^N\hookrightarrow\Sp_q^N$; compare \cite[Figure~2(a)]{PS2022}. Light grey indicates parameter regions with previously known sharp asymptotics. Dark grey highlights $0<p<1$ and $2<q\le\infty$, where Theorem~\ref{thm:main} closes the remaining gap and Remark~\ref{rem:all-indices} gives the resulting all-index asymptotics. The hatched triangle corresponds to $2<p<q\le\infty$, where matching bounds are still missing for intermediate indices $n$. The boundary of the two highlighted regions is not part of the open cases: sharp asymptotic estimates were already known for $p=1$ with $2<q\le\infty$, for $q=2$ with $0<p<1$, for $p=2$ with $2<q\le\infty$, and along the diagonal $p=q>2$.}
\label{fig:gelfand-regions}
\end{figure}

\begin{rem}[Complete asymptotic description]\label{rem:all-indices} 
Combining Theorem~\ref{thm:main} with \cite[Proposition~4.5]{HPV2021} and \cite[Proposition~3.1]{PS2022} yields a complete description of the Gelfand numbers for $0<p\le1$ and $2\le q\le\infty$. More precisely, with the universal constant $c_0\in(0,1)$ introduced above, we have 
\[
    c_n(\embed pq)
    \asymp_{p,q}
    \begin{cases}
        \displaystyle
        \min\left\{1,\frac{N}{n}\right\}^{1/p-1/2}
        &:\, 1\le n\le(1-c_0)N^2,\\[3mm]
        \displaystyle
        \max\left\{
            N^{1/q-1/p},
            N^{-1/p-1/2}\sqrt{N^2-n+1}
        \right\}
        &:\, (1-c_0)N^2<n\le N^2.
    \end{cases}
\]
The estimates in the second range were already sharp. The contribution of Theorem~\ref{thm:main} is the missing matching lower bound in the first range for $0<p<1$ and $q>2$. Moreover, in this first range, both comparison constants can be chosen independently of $q$.
\end{rem}

The remaining gap for Gelfand numbers identified in \cite{PS2022} remains open. In the notation used here, \cite[Conjecture~1.3]{PS2022} predicts that, for $2<p<q\le\infty$ and
\[
    (1-c_0)N^2
    \le n
    \le
    N^2-c_0N^{1+2/q}+1,
\]
one has
\[
    c_n(\embed pq)
    \lesssim_{p,q} \sqrt{\frac{N^2-n+1}{N^2}}^{\frac{1/p-1/q}{1/2-1/q}}. 
\]
Note that the original conjecture in \cite{PS2022} is slightly stronger: it asserts the same upper bound for the corresponding approximation numbers.
 
\subsubsection*{Idea of the Proof of Theorem~\ref{thm:main}}

The starting point is the rank-one case. If $n\le N$, an arbitrary subspace of $\R^{N\times N}$ of codimension less than $n$ must intersect nontrivially the $N$-dimensional space of matrices supported on a single column. Since every nonzero matrix in this space has rank one, all Schatten quasi-norms coincide there. For larger codimensions, this space is too small. It is therefore natural to replace it by the space of matrices supported on $r$ columns, which we identify with $\R^{N\times r}$. This space has dimension $Nr$, while all its elements have rank at most $r$. Thus, in order to force a nontrivial intersection with a subspace of codimension less than $n$, one is led to choose $Nr$ of order $n$, or equivalently $r$ of order $n/N$. This already explains the scale entering the desired estimate.

A rank bound alone, however, would lose a factor of order $\sqrt r$. The key observation is that the Hilbert--Schmidt geometry of the intersection produces a particularly well-positioned matrix. Indeed, projecting the standard basis matrices of $\R^{N\times r}$ onto a $d$-dimensional subspace and averaging their squared Hilbert--Schmidt norms yields one projected matrix $A$ for which
\[
    \|A\|_\infty\ge \|A\|_2^2.
\]
Since every matrix in $\R^{N\times r}$ has rank at most $r$, we have $\|A\|_p\le r^{1/p-1/2}\|A\|_2$. Together with $\|A\|_\infty\ge \|A\|_2^2$, this yields a lower bound of order $r^{1/2-1/p}$ whenever the subspace occupies a fixed proportion of $\R^{N\times r}$. Choosing $r\asymp n/N$ then produces precisely the factor
\[
    \left(\frac{N}{n}\right)^{1/p-1/2}.
\]
Finally, the passage from the operator norm to all target exponents $q\ge2$ uses $\|A\|_q\ge\|A\|_\infty$.

\subsubsection*{Organization of the paper}

The paper is organized as follows. Section~\ref{sec:preliminaries} contains the notation and elementary facts about Schatten norms and Gelfand numbers. In Section~\ref{sec:projection}, we establish the rectangular projection estimate and prove Theorem~\ref{thm:main}. 

\section{Notation and preliminaries}\label{sec:preliminaries}

We present some notation and preliminary facts that are used in the paper.

\paragraph{Matrix spaces and orthogonal projections.} 
For positive integers $a,b$, write $\M_{a,b}:=\R^{a\times b}$ and $\M_N:=\M_{N,N}$, and $\Tr$ for the trace of a linear operator on a finite-dimensional vector space. We equip $\M_{a,b}$ with the Hilbert--Schmidt inner product
\[
    \langle A,B\rangle_{\HS}:=\Tr(A^{\mathsf T}B),
    \qquad A,B\in\M_{a,b}.
\]
Let $e_1,\ldots,e_a$ and $f_1,\ldots,f_b$ denote the standard basis vectors of $\R^a$ and $\R^b$, respectively, regarded as column vectors. The matrices
\[
    E_{ij}:=e_i f_j^{\mathsf T},
    \qquad 1\le i\le a,\quad 1\le j\le b,
\]
each having a single entry equal to one in position $(i,j)$ and zeros elsewhere, form an orthonormal basis of $\M_{a,b}$. For a linear subspace $F\subset\M_{a,b}$, we write $P_F$ for the orthogonal projection onto $F$ with respect to this inner product. Viewed as a linear operator on $\M_{a,b}$, it has trace $\dim F$.

\paragraph{Schatten quasi-norms and rank estimates.}
For $A\in\M_{a,b}$, let $\sigma_1(A)\ge\cdots\ge\sigma_{\min\{a,b\}}(A)\ge0$ denote its singular values. We extend the notation $\|A\|_p$ from square to rectangular matrices by taking the $\ell_p$-norm, or quasi-norm when $0<p<1$, of this sequence. In particular, $\|A\|_\infty=\sigma_1(A)$ is the operator norm of $A:\ell_2^b\to\ell_2^a$, and $\|A\|_2^2=\langle A,A\rangle_{\HS}$. 

If $\rank A\le r$, where $r\ge1$, then for $0<p\le q\le\infty$
\begin{equation}\label{eq:rank-comparison}
    \|A\|_q
    \le\|A\|_p
    \le r^{1/p-1/q}\|A\|_q.
\end{equation}
The first inequality follows from the monotonicity of sequence quasi-norms. For the second, note that there are at most $r$ nonzero singular values. If $0<p<q<\infty$, H\"older's inequality with exponents $q/p$ and $q/(q-p)$ gives, after padding the nonzero singular-value sequence with zeros,
\[
    \sum_{j=1}^{r}\sigma_j(A)^p
    \le
    \left(\sum_{j=1}^{r}\sigma_j(A)^q\right)^{p/q}
    r^{1-p/q}.
\]
Taking the $p$-th root proves the assertion. The cases $p=q$ and $q=\infty$ are immediate. Under the same rank assumption, we shall in particular use that for all $0<p\le2$
\begin{equation}\label{eq:rank-p2}
    \|A\|_p
    \le r^{1/p-1/2}\|A\|_2
\end{equation}

\paragraph{Rectangular matrices as subspaces.}
For $1\le r\le N$, we identify $\M_{N,r}$ with the subspace of $\M_N$ consisting of matrices whose last $N-r$ columns vanish, via $A\mapsto\widetilde A:=(A\;0)$. Since $\widetilde A\widetilde A^{\mathsf T}=AA^{\mathsf T}$, this identification leaves all nonzero singular values unchanged. Consequently, $\|\widetilde A\|_p=\|A\|_p$, $0<p\le\infty$. We shall therefore pass between these rectangular and square matrix representations without changing the notation.

\paragraph{Gelfand numbers.}
For $1\le n\le N^2$, the condition $\codim F<n$ in \eqref{eq:def-gelfand} is equivalent to $\dim F\ge N^2-n+1$. Restricting to a smaller subspace cannot increase the restriction norm. Hence the infimum may be taken over subspaces of dimension exactly $N^2-n+1$, and
\begin{equation}\label{eq:dimension-form}
    c_n(\embed pq)
    =
    \inf_{\substack{F\subset\M_N\\
                     \dim F=N^2-n+1}}
    \sup_{A\in F\setminus\{0\}}
    \frac{\|A\|_q}{\|A\|_p}.
\end{equation}

For $0<p\le q\le\infty$, the full embedding has norm one: the upper bound follows from \eqref{eq:rank-comparison}, and equality is attained on any nonzero rank-one matrix. Moreover, applying \eqref{eq:rank-comparison} with $r=N$ gives
\begin{equation}\label{eq:dimension-floor}
    c_n(\embed pq)\ge N^{1/q-1/p},
    \qquad 1\le n\le N^2.
\end{equation}
This lower bound is attained at the last index. Indeed, let $I_N=\operatorname{diag}(1,\ldots,1)$ denote the $N\times N$ identity matrix and choose $F=\operatorname{span}\{I_N\}$ in \eqref{eq:dimension-form} for $n=N^2$. All singular values of $tI_N$ equal $|t|$, so for every $t\ne0$,
\[
    \frac{\|tI_N\|_q}{\|tI_N\|_p}
    =N^{1/q-1/p}.
\]
Together with \eqref{eq:dimension-floor}, for $0<p\le q\le\infty$ this yields
\[
    c_{N^2}(\embed pq)=N^{1/q-1/p}.
\]
Finally, the inequality $\|A\|_q\ge\|A\|_\infty$ implies
\[
    c_n(\embed{p}{q})
    \ge c_n(\embed{p}{\infty}),
    \qquad 0<p,q\le\infty,\quad 1\le n\le N^2.
\]
Thus, to prove the lower bound in Theorem~\ref{thm:main}, it suffices to consider the target space $\Sp_\infty^N$.

\section{A rectangular projection argument}\label{sec:projection}

We first establish the subspace estimate that is the main ingredient of the proof. Its formulation includes the range $0<p\le2$, although only $0<p\le1$ is needed for Theorem~\ref{thm:main}.

\begin{lemma}\label{lem:rectangular}
Let $1\le r\le N$, and let $E\subset\M_{N,r}$ be a subspace of dimension $d\ge1$. Then, for $0<p\le2$,
\begin{equation}\label{eq:rectangular}
 \sup_{A\in E\setminus\{0\}}
 \frac{\|A\|_\infty}{\|A\|_p}
 \ge\frac{\sqrt d}{\sqrt N\,r^{1/p}}.
\end{equation}
\end{lemma}

\begin{proof}
Let $P_E:\M_{N,r}\to\M_{N,r}$ be the
Hilbert--Schmidt-orthogonal projection onto the subspace $E$. Since $P_E=P_E^{*}=P_E^2$, we have
\[
    \|P_EE_{ij}\|_2^2
    =
    \langle P_EE_{ij},E_{ij}\rangle_{\HS}.
\]
As the trace of a finite-dimensional linear operator is independent of the choice of basis and the matrices \(E_{ij}\) form an orthonormal basis of \(\M_{N,r}\), we have
\[
    \Tr(P_E)
    =
    \sum_{i=1}^N\sum_{j=1}^r
    \langle P_EE_{ij},E_{ij}\rangle_{\HS}.
\]
But since $P_E$ is the orthogonal projection onto the $d$-dimensional space $E$, we have $\Tr(P_E)=d$ and, consequently,
\[
    \sum_{i=1}^N\sum_{j=1}^r\|P_EE_{ij}\|_2^2=d.
\]
Hence, by averaging, there must exist indices $i,j$ such that
\[
    t:=\|P_EE_{ij}\|_2^2\ge\frac d{Nr}.
\]
Now put $A:=P_EE_{ij}=(a_{k\ell})_{k,\ell}$. Since \(P_E\) is the orthogonal projection onto \(E\), we have
\[
A=P_EE_{ij}\in E
\qquad\text{and}\qquad
E_{ij}-A = \bigl(I_{\M_{N,r}}-P_E\bigr)E_{ij} = P_{E^\perp} E_{ij}
\in E^\perp,
\]
and so $\langle A,E_{ij}-A\rangle_{\HS}=0$. Therefore
\[
    a_{ij}
    =\langle A,E_{ij}\rangle_{\HS}
    =\langle A,A\rangle_{\HS}
    =\|A\|_2^2
    =t.
\]
Since the absolute value of every matrix entry is bounded by the operator norm,
\[
    \|A\|_\infty\ge |a_{ij}|=t.
\]
On the other hand, \(A\in\M_{N,r}\) and hence \(\rank A\le r\). By \eqref{eq:rank-p2},
\[
    \|A\|_p
    \le r^{1/p-1/2}\|A\|_2
    =r^{1/p-1/2}\sqrt t.
\]
It follows that
\[
    \frac{\|A\|_\infty}{\|A\|_p}
    \ge r^{1/2-1/p}\sqrt t
    \ge r^{1/2-1/p}\sqrt{\frac d{Nr}}
    =\frac{\sqrt d}{\sqrt N\,r^{1/p}},
\]
which proves the result.
\end{proof}

The parameter $r$ determines the dimension $Nr$ of the rectangular matrix space. Intersecting it with a subspace of codimension less than $n$ leaves a subspace of dimension at least $Nr-n+1$. We keep $r$ free because the resulting lower bound depends on $r$, and a suitable choice of $r$ later on will balance the intersection dimension against the rank-dependent norm estimate.

\begin{proposition}\label{prop:all-index-lower} 
Let $0<p\le2\le q\le\infty$, $N\in\N$, and $1\le n\le N^2$. Then
\begin{equation}\label{eq:all-index-lower}
 c_n(\embed pq)
 \ge
 \max_{\substack{1\le r\le N\\Nr\ge n}}
 \frac{\sqrt{Nr-n+1}}{\sqrt N\,r^{1/p}}.
\end{equation}
Moreover, $c_n(\embed pq)=1$ whenever $1\le n\le N$.
\end{proposition}

\begin{proof}
Fix a subspace $F\subset\M_N$ with $\codim F<n$. For $1\le r\le N$, let $H_r\subset\M_N$ denote the subspace of matrices whose last $N-r$ columns vanish. By deleting these zero columns, we identify $H_r$ with $\M_{N,r}$. As observed above, this identification preserves all Schatten norms and quasi-norms, and $\dim H_r=Nr$. Now consider the canonical quotient map
\[
    Q:\M_N\longrightarrow \M_N/F,\qquad Q(A)=A+F.
\]
Its restriction \(Q|_{H_r}:H_r\to\M_N/F\) has kernel $\ker(Q|_{H_r})=F\cap H_r$. Moreover, $\rank(Q|_{H_r}) \le \dim(\M_N/F) =\codim F$ and
hence, by the rank--nullity theorem, 
\begin{equation}\label{eq:intersection-dimension}
    \dim(F\cap H_r)
    =\dim H_r-\rank(Q|_{H_r})
    \ge Nr-\codim F
    \ge Nr-n+1.
\end{equation}

Suppose now that $Nr\ge n$. Then $F\cap H_r\ne\{0\}$. Put $d:=\dim(F\cap H_r)$. Applying Lemma~\ref{lem:rectangular} to the subspace $F\cap H_r\subset H_r\cong\M_{N,r}$, and using $d\ge Nr-n+1$, we obtain
\[
    \sup_{A\in(F\cap H_r)\setminus\{0\}}
    \frac{\|A\|_\infty}{\|A\|_p}
    \ge
    \frac{\sqrt d}{\sqrt N\,r^{1/p}}
    \ge
    \frac{\sqrt{Nr-n+1}}{\sqrt N\,r^{1/p}}.
\]
Since $F\cap H_r\subset F$ and $\|A\|_q\ge\|A\|_\infty$ for every $q>0$, it follows that
\[
    \sup_{A\in F\setminus\{0\}}
    \frac{\|A\|_q}{\|A\|_p} 
    \ge
    \frac{\sqrt{Nr-n+1}}{\sqrt N\,r^{1/p}}.
\]
This estimate holds for every integer $r$ satisfying $1\le r\le N$ and $Nr\ge n$. We may therefore maximize over all such $r$. Since $F$ was arbitrary, taking the infimum over all subspaces $F\subset\M_N$ with $\codim F<n$ yields
\[
    c_n(\embed pq)
    \ge
    \max_{\substack{1\le r\le N\\ Nr\ge n}}
    \frac{\sqrt{Nr-n+1}}{\sqrt N\,r^{1/p}},
\]
which proves \eqref{eq:all-index-lower}.

Finally, suppose that $1\le n\le N$. Taking $r=1$ in \eqref{eq:intersection-dimension} gives $\dim(F\cap H_1)\ge N-n+1\ge1$, and hence $F\cap H_1\ne\{0\}$. Every nonzero matrix in $H_1$ has rank one, so all of its Schatten norms and quasi-norms coincide. Therefore
\[
    \sup_{A\in F\setminus\{0\}}
    \frac{\|A\|_q}{\|A\|_p}\ge1.
\]
Since $p\le q$ in the parameter range under consideration, the full embedding $\Sp_p^N\hookrightarrow\Sp_q^N$ has norm one. Consequently, whenever $1\le n\le N$,
\[
    c_n(\embed pq)=1,
\]
which completes the proof.
\end{proof}

We now choose the rectangular dimension in Proposition~\ref{prop:all-index-lower} and prove the main result.

\begin{proof}[Proof of Theorem~\ref{thm:main}]
Let us define $\alpha:=1/p-1/2$. We first prove the lower estimate. The proof uses only $0<p\le2$, and we distinguish three ranges of indices.

\medskip
\noindent\textit{Case 1: $1\le n\le N$.}
Proposition~\ref{prop:all-index-lower} gives $c_n(\embed pq)=1$. Since $\min\{1,N/n\}=1$, the lower estimate holds with constant one.

\medskip
\noindent\textit{Case 2: $N<n\le N^2/2$.}
Here we choose $r:=\left\lceil\frac{2n}{N}\right\rceil$. Since $n\le N^2/2$, we have $r\le N$. Moreover, $Nr\ge2n$, and therefore
\[
    Nr-n+1\ge\frac{Nr}{2}.
\]
On the other hand,
\[
    r\le\frac{2n}{N}+1\le\frac{3n}{N},
\]
where the last inequality uses $n>N$. With our specific choice of $r$, we can apply the lower estimate in \eqref{eq:all-index-lower} combined with the previous estimates to obtain
\[ 
 c_n(\embed pq)
 \ge\frac{\sqrt{Nr/2}}{\sqrt N\,r^{1/p}}
 =\frac1{\sqrt2}r^{-\alpha}
 \ge\frac{3^{-\alpha}}{\sqrt2}
 \left(\frac Nn\right)^{\alpha}.
\]
In the last estimate we used $r\le \frac{3n}{N}$.

\medskip
\noindent\textit{Case 3: $N^2/2<n\le(1-\eta)N^2$.}
If this range is non-empty, take $r=N$ in \eqref{eq:all-index-lower}. Since in our range of $n$, $N^2-n+1\ge\eta N^2$, we obtain
\[
 c_n(\embed pq)
 \ge N^{-1/2-1/p}\sqrt{N^2-n+1}
 \ge\sqrt\eta\,N^{-\alpha}.
\]
The inequality $n>N^2/2$ implies
\[
 \min\left\{1,\frac Nn\right\}^{\alpha} 
 \le\left(\frac{2}{N}\right)^{\alpha}. 
\]
Consequently,
\[
 c_n(\embed pq)
 \ge\sqrt\eta\,2^{-\alpha}
 \min\left\{1,\frac Nn\right\}^{\alpha};
\]
empty ranges are omitted. 
\medskip

Combining the three cases above proves the lower estimate in \eqref{eq:main}; one may take 
\begin{equation}\label{eq:explicit-constant}
 c_{p,\eta}
 :=\min\left\{1,\frac{3^{-\alpha}}{\sqrt2},\sqrt\eta\,2^{-\alpha}\right\}.
\end{equation}
In particular, this constant is independent of $q$.

\medskip 
For the upper estimate, we use the known bound
\begin{equation}\label{eq:upper-hilbert}
 c_n(\embed p2)
 \le C_p\min\left\{1,\frac Nn\right\}^{\alpha},
 \qquad 1\le n\le N^2;
\end{equation}
see \cite[Proposition~3.3, with $q=2$]{PS2022}; compare also \cite[Proposition~4.5, Step~1]{HPV2021}. Since $q\ge2$, the inequality $\|A\|_q\le\|A\|_2$ holds for every matrix $A$. Taking the restriction norm on each admissible subspace and then the infimum gives
\[
 c_n(\embed pq)\le c_n(\embed p2).
\]
Together with \eqref{eq:upper-hilbert}, this proves the upper estimate in \eqref{eq:main} with a constant independent of both $q$ and $\eta$.
\end{proof}

\section*{Acknowledgement}
The author was supported by the German Research Foundation under DFG project 516672205 \emph{Limit theorems for the volume of random projections of $\ell_p$-balls}. Most of this work was carried out during a research semester at Graz University of Technology. The author gratefully acknowledges the excellent working conditions and the hospitality provided during his stay.
 
ChatGPT (GPT-5.6 Luna, OpenAI) was used to improve the wording and for grammar and spelling checks, and to assist in the production of Figure \ref{fig:gelfand-regions}.

\bibliographystyle{abbrv}
\bibliography{gelfand}

\end{document}